\documentclass{wart}
\usepackage{xspace,amsmath,palatino,amsthm,amsfonts,amssymb}
\newtheorem{thm}{Theorem}[section]
  \newtheorem*{thm*}{Theorem}
  \newtheorem{prop}[thm]{Proposition}
  
  \newtheorem{cor}[thm]{Corollary}
  \newtheorem{conj}{Conjecture}
  \newtheorem{question}{Question}
 
 \theoremstyle{remark}
 
 \newtheorem{remark}{Remark}

 \newcommand{\nc}{\newcommand}
 \newcommand{\renc}{\renewcommand}
 \nc{\daloc}{\tilde{\aloc}}
 \nc{\ex}{\mathbf{e}}
 \nc{\ASt}{\mathcal{AS}}
 \nc{\ASto}{\ASt_0}
 \nc{\gN}{N^V}
 \nc{\Bp}{B^+}
 \nc{\Bm}{B^-}
 \nc{\clus}[3]{E^{#1}_{#2,#3}}
 \nc{\drse}{double reduced subexpression\xspace}
 \nc{\ddse}{double distinguished subexpression\xspace}
 \nc{\ddses}{double distinguished subexpressions\xspace}
 \nc{\Gsp}{G^{>0}}
 \nc{\Gp}{G^{\geq 0}}
 \nc{\pos}[1]{{#1}^{\geq 0}}
 \nc{\rel}[1]{\overset{#1}\longrightarrow}
 \nc{\brel}[1]{\overset{#1}\longleftarrow}
 \nc{\Qi}[3]{\mathcal{Q}^{#1}_{#2,#3}}
 \nc{\Qb}[3]{\mathcal{Q}^{#1}_{\mathbf{#2,#3}}}
 \nc{\gdb}[3]{\mathcal{P}^{#1}_{#2,#3}}
 \nc{\gdbb}[3]{\mathcal{P}^{\mathbf{#1}}_{\mathbf{#2,#3}}}
 \nc{\rep}{representation}
 \nc{\un}{u_{(n)}}
 \nc{\wn}{w_{(n)}}
 \nc{\vn}{v_{(n)}}
 \nc{\tuk}{\tilde{u}_{(k)}}
 \nc{\tukm}{\tilde{u}_{(k-1)}}
 \nc{\unm}{u_{(n-1)}}
 \nc{\vnm}{v_{(n-1)}}
 \nc{\wnm}{w_{(n-1)}}
 \nc{\ukm}{u_{(k-1)}}
 \nc{\vkm}{v_{(k-1)}}
 \nc{\wkm}{w_{(k-1)}}
 \nc{\uwkm}{\uw_{(k-1)}}
 \nc{\uvkm}{\uv_{(k-1)}}
 \nc{\uwk}{\uw_{(k)}}
 \nc{\uvk}{\uv_{(k)}}
 \nc{\uw}{u^\Bw}
 \nc{\uv}{u^\Bv}
 \nc{\Bt}{\mathbf{t}}
 \nc{\bw}{\bar w}
 \nc{\bu}{\bar u}
 \nc{\bv}{\bar v}
 \nc{\Bk}{B^{(k)}}
 \nc{\won}{\bar{G}}
 \nc{\xvw}{x_{\mathbf{v,w}}}
 \nc{\CsC}[2]{(\C^*)^{#1}\times\C^{#2}}
 \nc{\hw}{highest weight\xspace} 
 \renc{\le}{\left}
 \nc{\ri}{\right}
 \nc{\al}{\alpha}
 \nc{\Bv}{\mathbf{v}}
 \nc{\Bw}{\mathbf{w}}
 \nc{\C}{\mathbb{C}}
 \nc{\De}{\Delta}
 \nc{\ep}{\epsilon}
 \nc{\fr}[1]{\mathfrak{#1}}
 \nc{\gm}[3]{\Delta^{#1\omega_{#3}}_{#2\omega_{#3}}}
 \nc{\inn}[1]{\langle #1\rangle}
 \nc{\Jo}[1]{J^{\circ}_{\mathbf{#1}}}
 \nc{\Jm}[1]{J^{-}_{\mathbf{#1}}}
 \nc{\Jp}[1]{J^{+}_{\mathbf{#1}}}
 \nc{\la}{\lambda}
 \nc{\mc}[1]{\mathcal{#1}}
 \nc{\Ob}{\mathrm{Ob}}
 \nc{\Om}{\Omega}
 \nc{\om}{\omega}
 \nc{\Orb}{\mathrm{Orb}}
 \nc{\pad}{\hat{\Z}_p}
 \nc{\pder}[2]{\frac{\partial #1}{\partial #2}}
 \nc{\pderw}[1]{\frac{\partial}{\partial #1}}
 \nc{\pdersec}[2]{\frac{\partial^2 #1}{\partial {#2}^2}}
 \nc{\perm}[1]{\pi_{#1}}
 \nc{\Q}{\mathbb{Q}}
 \nc{\qvw}[1]{\La(#1 \Bv,\Bw)}
 \nc{\qK}{\frac{K-K^{-1}}{q-q^{-1}}}
 \nc{\R}{\mathbb{R}}
 \nc{\rad}{\mathrm{rad}\,}
 \nc{\Rep}{\mathrm{Rep}}
 \nc{\res}[2]{\mathrm{res}^{#1}_{#2}}
 \nc{\resp}[2]{\mathfrak{res}^{#1}_{#2}}
 \nc{\rg}{\Sigma^1}
 \nc{\RG}{\EuScript{R}_G}
 \nc{\Ri}[2]{\mathcal{R}_{#1,#2}}
 \nc{\Rb}[2]{\Ri{\mathbf{#1}}{\mathbf{#2}}}
 \nc{\HH}{\mathbb{H}}
 \nc{\rank}{\mathrm{rk}\,}
 \nc{\rmat}{$R$-matrix\xspace}
 \nc{\V}[1]{\mathbf{#1}}
 \nc{\vp}{\varphi}
 \nc{\scm}[1]{\Delta^{v_{(#1)}\om_{i_{#1}}}_{w_{(#1)}\om_{i_{#1}}}}
 \nc{\spcm}[1]{\Delta^{v_{(#1-1)}\om_{i_{#1}}}_{w_{(#1)}\om_{i_{#1}}}}
 \nc{\si}{\sigma}
 \nc{\sk}{\mathrm{sk}\,}
 \nc{\Stab}{\mathrm{Stab}}
 \nc{\SL}[1]{\mathrm{SL}_{#1}}
 \nc{\SO}[1]{\mathrm{SO}_{#1}}
 \nc{\so}[1]{\fr{so}_{#1}}
 \nc{\Sp}[1]{\mathrm{Sp}_{#1}}
 
 \nc{\Spec}{\mathrm{Spec}}
 \nc{\SU}[1]{\mathrm{SU}( #1)}
 \nc{\su}[1]{\fr{su}_{#1}}
 \nc{\Sym}{\mathrm{Sym}}
 \nc{\sym}{\mathrm{sym}}
 \nc{\Tg}{\mc{T}(\fr g)}
 \renc{\th}{\theta}
 \nc{\ttH}{{_2 H}}
 \nc{\thi}{\th^{-1}}
 \nc{\tG}{\tilde{G}}
 \nc{\tom}{\tilde{\omega}}
 \nc{\tpi}{\tilde{\pi}}
 \nc{\Tr}{\mathrm{Tr}}
 \nc{\tr}{\mathrm{tr}}
 \nc{\tD}{\tilde{\Delta}}
 \nc{\tri}{\tau}
 \nc{\Ug}{\mc{U}(\fr g)}
 \nc{\Un}{\mc{U}(\fr n)}
 \nc{\Uh}{\mc{U}(\fr h)}
 \nc{\wt}{\mathrm{wt}}
 \renewcommand{\V}[1]{\mathbf{#1}}
 \nc{\Z}{\mathbb{Z}}
 \nc{\Zp}{\Z/p}
 \nc{\Znn}{\Z_{\geq 0}}

\theoremstyle{remark}

\nc{\flag}{\mathcal{B}}
\nc{\vk}{v_{(k)}}
\nc{\wk}{w_{(k)}}
\nc{\uk}{u_{(k)}}
\nc{\Bu}{\mathbf{u}}
\nc{\Up}{U^+}
\nc{\Um}{U^-}
\begin{document}
\title[A Deodhar-type stratification]{A
  Deodhar-type stratification\\  on the double flag variety}
\authors{\textsc{Ben Webster}&\textsc{Milen Yakimov}}
\address{Department of Mathematics &Department of Mathematics\\
         University of California & University of California\\
          Berkeley, CA 94720 &Santa Barbara, CA 93106}
\email{bwebste@math.berkeley.edu}
\urladdr{\tt{http://math.berkeley.edu/~bwebste}}

\email{yakimov@math.ucsb.edu}


\begin{abstract}
  We describe a partition of the double flag variety 
  $G/\Bp\times G/\Bm$ of a complex semisimple algebraic group $G$
  analogous to the Deodhar partition on the flag variety
  $G/\Bp$. This partition is a refinement of the stratification into orbits both
  for $\Bp\times\Bm$ and for the diagonal action of $G$, just as
  Deodhar's partition refines the orbits of $\Bp$ and $\Bm$.
  
  We give a coordinate system on each stratum, and show that all strata
  are coisotropic subvarieties.  
  Also, we discuss possible connections to the positive and cluster
  geometry of $G/ \Bp \times G/ \Bm$, which would generalize results of
  Fomin and Zelevinsky on double Bruhat cells and Marsh and Rietsch on
  double Schubert cells.
\end{abstract}

\maketitle

\section{Introduction}
\label{sec:introduction}
Let $G$ be a complex semisimple algebraic group with a fixed pair of
opposite Borel subgroups $B^\pm$ and let $T=B^+ \cap B^-$ be the
corresponding maximal torus of $G$. The Weyl group of the pair $(G, T)$
will be denoted by $W$ and its identity element by $e$. The flag variety,
that is the set of Borel subgroups of $G$, is naturally identified with $G/B$ for
any Borel subgroup $B$ by the map $gB\mapsto gBg^{-1}$ and naturally has the
structure of a projective variety.

In \cite{Deo85}, Deodhar described a remarkable partition of the
flag variety $G/\Bp=\flag$.  It is a refinement of the stratification of
$\flag$ both into $\Bp$ and into $\Bm$ orbits, and thus of the
stratification of $\flag$ into intersections of 
opposite Schubert cells $\Ri vw=\Bp w\cdot\Bp\cap \Bm v\cdot\Bp$
which we call {\em{double Schubert cells}}. In particular, it provides a
stratification of each double Schubert cell.
Furthermore, he showed that each stratum is isomorphic to $\C^k\times
(\C^*)^m$, where $m+2k=\dim \Ri vw=\ell(w)-\ell(v)$, and the strata are
indexed by certain subexpressions of a reduced expression for $w$.

Throughout the paper by a {\em{stratification}} of a quasiprojective
variety $X$ we mean a partition of $X$
\[
X = \bigsqcup_{\alpha \in A} X_{\alpha} 
\]
into smooth, locally closed, irreducible subsets $X_\alpha$, 
the closure of each of which is a union of strata $X_\beta$,
$\beta \in A(\alpha) \subset A$.  

In this paper, we will describe an analogue of Deodhar's partition
in the double flag variety $\flag\times\flag$.

Playing the role of the subgroups $\Bp$ and $\Bm$ are the actions of
$G_\De$ (the diagonal subgroup of $G \times G$)
and $\Bp\times \Bm$ as subgroups of $G\times G$.  These actions
are in a certain sense (Poisson) dual, much as $\Bp$ and $\Bm$ are.
To be more precise, recall that $G \times G$ is a double Poisson algebraic 
group for the standard Poisson structure on $G$. The image of $G$ inside 
$G \times G$ is $G_\De$ and the normalizer of its dual is $\Bp \times \Bm$.
Similarly a double of $\Bp$ is closely related to $G$ and its dual is 
$\Bm$. 

As is well known, the orbits of $G_\De$ on $\flag \times \flag$
are indexed by the Weyl group $W$ of
$G$, and $\Bp\times\Bm$-orbits are indexed by $W\times W$, both with the
partial order induced by closure being the standard strong Bruhat order
(or its opposite, depending on indexing conventions).
Thus, to each ordered triplet $(u,v,w)\in W^3$, we associate the
intersection of the $G_\De$ orbit attached to $u$ and
$\Bp\times\Bm$-orbit associated to $(v,w)$.
\begin{equation}
  \gdb uvw=\le(G_\De\cdot (\Bp,u\cdot \Bm)\ri)\bigcap\le(\le(\Bp\times \Bm\ri)\cdot
  (v\cdot\Bp,w\cdot\Bm)\ri).  
\label{0:part}
\end{equation}

We will prove a number of results about these remarkable varieties.
Here, we collect the principal results of Theorem \ref{thm:strata-param}
and Proposition \ref{prop:coisotropy}:
\begin{thm*}
  For all $u,v,w\in W^3$, the variety $\gdb uvw$ has a stratification
  indexed by combinatorial objects related to $u,v,w$ called \ddses,
  such that each stratum has explicit coordinates, which realize an
  isomorphism with $\CsC{k}{\ell}$ for some integers $k,\ell$.
  Furthermore, each stratum is coisotropic in the double Poisson
  structure on $\flag\times\flag$. 
\end{thm*}

Interest in these varieties has arisen from different directions in
mathematics. In \cite[Example 4.9]{EL04}, they were identified with the 
torus orbits of symplectic leaves of $\flag\times\flag$
considered as a Poisson subvariety of variety of Lagrangian
subalgebras associated to the standard Poisson structure on $G$.

On the other hand, similar varieties appear in the study of positivity
in algebraic groups and cluster algebras.  The orbit intersection $\gdb
{e}vw$ is just a reduced double Bruhat cell as defined by Fomin and
Zelevinsky in \cite{FZ99}, and thus has a natural cluster algebra
structure.  Furthermore, the intersections with double Bruhat cells give
a cell decomposition of the non-negative part of $G$. 
The stratum $\gdb {w_0}vw$ is naturally a double Schubert cell, and the
intersections of the non-negative part of $\flag$ with double Schubert
cells gives a similar cell decomposition.

When $G$ is of adjoint type, the double flag variety $\flag \times
\flag$ sits inside the wonderful compactification $\overline{G}$ of $G$
as the so called lowest stratum (the unique closed $G \times G$ orbit).
For the purposes of understanding total positivity in $\overline{G}$
better, Lusztig \cite{Lus98b} defined a partition of $\overline{G}$
which is a refinement of of the partition by $G \times G$ orbits. The
partition \eqref{0:part} differs from the restriction of Lusztig's
partition to $\flag \times \flag$. The latter consists of products of
Schubert cells.

In Section~\ref{sec:relat-posit-flagt}, we will define the varieties in
question, and discuss their relationship to double Schubert cells and
double Bruhat cells.  In Section~\ref{sec:bruhat-decomposition}, we will
recall the background on relative positions and Bruhat decomposition
that we will need later, and in Section~\ref{sec:ddses}, we cover the
analogue, for our varieties, of Deodhar's theory of distinguished
subexpressions, which provide the necessary combinatorial setup.  In
Section~\ref{sec:strat-gdb-uvw}, we will define our stratification and
define a parameterization of each stratum of it, as well as
characterization of each stratum in terms of generalized minors in
Section~\ref{sec:chamber-minors}.  In
Sections~\ref{sec:positivity-1}~and~\ref{sec:coisotr-deodh-strat}, we
will consider connections to positivity and Poisson geometry respectively.

\section*{Acknowledgments}
\label{sec:acknowledgements}

We would like to thank Nicolai Reshetikhin
for his help and encouragement. We are grateful to the Managing 
Editors for bringing the very interesting paper \cite{Cur88} to our 
attention.

B.W. was supported by a National Science Foundation 
Graduate Research Fellowship and by the RTG grant DMS-0354321.

M.Y. was supported by the National
  Science Foundation grant DMS-0406057 and an Alfred P. Sloan research
  fellowship.

\section{The varieties $\gdb uvw$}
\label{sec:relat-posit-flagt}

The varieties $\gdb uvw$ can be defined in terms of the more classical
terminology of relative position.  For any two flags $B',B''\in\flag$,
we can define a ``distance'' (in a very crude sense) between them as
follows: The orbits of $B'$ on $\flag$ are naturally indexed by $W$, and
relative position $r(B',B'')\in W$ indexes the orbit in which $B''$
lies.  Following \cite{MR04}, we use the notation $B\rel{w}B'$ to
express that $B$ and $B'$ are in relative position $w$.  Alternatively,
relative positions can be characterized as the $G_\De$-orbits on
$\flag\times\flag$. That is, $B\rel{w} B'$ if and only if $(g\cdot
B,g\cdot B')=(\Bp ,w\cdot \Bp )$ for some $g\in G$.  Relative positions
can be used to define several interesting subvarieties of $\flag$ and
$\flag\times\flag$.

For ease of notation let $v^*=w_0vw_0$, where $w_0\in W$ is the longest
element of $W$.

\begin{itemize}
\item The Schubert cells of
$\flag$ with respect to $\Bp$ and $\Bm$ can be defined by:
\begin{align*}
  \flag_w&=\Bp w\cdot \Bp=\{B\in\flag \vert \Bp \rel{w} B\},\\
  \flag^v&=\Bm v\cdot \Bm=\{B\in\flag \vert \Bm \rel{v^*} B\} 
\end{align*}
\item The double Schubert cells can be defined by
\begin{equation}\label{eq:2}
  \Ri wv=\flag_v\cap\flag^{ww_0}=\{B\in\flag \vert \Bp \rel{v} B \brel{w_0w} B^-\}.
\end{equation}

\item In $\flag\times\flag$, as we described above 
\begin{align*}
  G_\De\cdot(\Bp,u\cdot \Bm)&=\{(B_1,B_2)\in\flag\times\flag\vert B_1
  \rel{uw_0} B_2\}\\
  (\Bp\times\Bm)\cdot(v\cdot\Bp,w\cdot\Bm)&=\{(B_1,B_2)\in\flag\times\flag\vert
  \Bp \rel{v}B_1,\Bm  
  \rel{w^*}  B_2\}
\end{align*}

\item Thus, the intersection of these orbits is the variety
\begin{equation}\label{eq:1}
  \gdb uvw=\{(B_1,B_2)\in\flag\times\flag\vert \Bp \rel{v}B_1 \rel{uw_0} B_2
  \brel{w^*} \Bm \}
\end{equation}
\end{itemize}

While this choice of indexing of relative positions may seem a little
strange, we will see in
Sections~\ref{sec:ddses}~and~\ref{sec:strat-gdb-uvw} that this will
simplify the combinatorics necessary to describe these varieties.

We consider some already known special cases:
For all $v,w\in W$, we let 
\begin{equation*}G_{v,w}=\Bp w \Bp\cap \Bm v \Bm \subset G\end{equation*} 
denote a double Bruhat cell
of $G$, and let $L_{v,w}=G_{v,w}/T$ be the corresponding reduced double
Bruhat cell.
\begin{prop}\label{double-schub}
  If $u=e$, then $\gdb {e}vw$ is
  isomorphic to $L_{v,w}$.  Similarly, the varieties $\gdb {w_0}v{w}\cong \gdb
  {w_0w^{-1}}ve\cong \gdb {vw_0}{e}{w}$ are isomorphic to $\Ri {ww_0}{v}$.
\end{prop}
\begin{proof}
  Rewriting the definition of a double Bruhat cell, any $g\in G$ is in
  $G_{u,v}$ if and only if $\Bp \rel{w}g\cdot \Bp $ and $\Bm \rel{v^*}
  g\cdot \Bm$.
  
  Since $g\cdot \Bp \rel{w_0}g\cdot \Bm $, the image of the action map
  sending $g\in G_{v,w}$ to $(g\cdot \Bp,g\cdot\Bm)$ is contained in
  $\gdb {e}vw$, and this map is surjective (since $G$ acts transitively
  on pairs of flags in fixed relative position).  Since $T$ is the
  stabilizer of $(\Bp,\Bm)$, this map is an isomorphism
  $L_{v,w}=G_{v,w}/T\cong \gdb {e}vw$.
  
  On the other hand, if either of $v,w$ is $e$, or $u=w_0$ then the
  definition of equation \eqref{eq:1} collapses to the definition in
  \eqref{eq:2}, for the appropriate choice of new $v,w$.
\qed\end{proof}

\section{Relative position and Bruhat decomposition}
\label{sec:bruhat-decomposition}

In this section, we will restate many well-known properties of relative
position and the Bruhat stratification of $\flag$, which will be used later.

\begin{prop}\label{prop:rel-pos-prop}
Relative position is 
\begin{enumerate}
\item invariant under the diagonal action of $G$,
  i.e. 
  \begin{equation*}
    r(B',B'')=r(g\cdot B',g\cdot B''),
  \end{equation*}
\item anti-symmetric, i.e. 
  \begin{equation*}
    r(B',B'')=\left( r(B'',B') \right)^{-1},
  \end{equation*}
\item sub-multiplicative with respect to Bruhat order, i.e. there
  exist $u,v\in W$ such that  $u\leq r(B',B''),v\leq r(B'',B''')$ and
  \begin{equation*}
    uv= r(B',B'''),
  \end{equation*}
\item\label{item:4} multiplicative if
  \begin{equation*}
    \ell\le(r(B',B'')\ri)+\ell\le(r(B'',B''')\ri) =
    \ell\le(r(B',B'')\cdot r(B'',B''') \ri).
  \end{equation*}
\end{enumerate}
\end{prop}

Complementary to property~\ref{item:4} is the following essential
result:

\begin{prop}
  If $B\rel{w}B'$, and $w=w'w''$ for $\ell(w')+\ell(w'')=\ell(w)$, then
  there is a unique flag $\tilde B$ such that $B\rel{w'}\tilde B\rel{w''}B'$.
\end{prop}

In the case where $B=\Bp$, we can use this to define a surjective
projection map $\pi^w_{w'}:\flag_w\to\flag_{w'}$. We will also use this
notation for analogous projections between $\Bp\times\Bm$-orbits on
$\flag\times\flag$. This makes $\flag_w$ a trivial fiber bundle over
$\flag_{w'}$.  In the case where $w's_i=w$, we can explicitly trivialize
this bundle.  This will be done in Proposition~\ref{Bruhat-para}. The
general case can be treated by iteration.

As usual, we let $x_{i}(t)=\exp(te_i)$ and
$x_{-i}(t)=y_i(t)=\exp(tf_i)$, where $e_i,f_i$ are the standard
Chevalley generators for the simple root spaces, and $s_i$ are the
corresponding simple reflections in the Weyl group.

Following Fomin and Zelevinsky, we pick a lift of $W$ to $G$ as follows:
for simple reflections, we let 
$\bar{s_i}=x_i(-1)y_i(1)x_i(-1)$, and $\bar{w}=\bar{w}_1\bar w_2$ if
$w=w_1w_2$ and \linebreak $\ell(w)=\ell(w_1)+\ell(w_2)$.

Note that this is {\bf not} a homomorphism, since $\bar s_i^2\neq 1$.
However, it is independent of the reduced word we choose.

We will use the notation $x_{w(\al_i)}(t)=\bar{w}x_i(t)\bar{w}^{-1}$.

\begin{prop}\label{Bruhat-para} Fix two flags 
$B_1$ and $B_2$ of relative position $B_1\rel{u} B_2$; that is
$(B_1,B_2)=g\cdot (u_1\cdot \Bp,u_2\cdot \Bp)$ for some $u_1,u_2 \in W$
such that $u=u_1^{-1}u_2$. Let $X\subset\flag$ be the
subvariety of flags $B'$ such that $B_1\rel{s_i}B'$.  Then $X\cong \C$,
and the following hold.
\begin{enumerate}
\item If $s_iu > u$, then $B'\rel{s_iu} B_2$ for all $B'\in X$ and the map 
\begin{equation*}
t\mapsto gx_{u_1(\al_i)}(t)u_1s_i\cdot \Bp
\end{equation*}
defines an isomorphism $\C\to X$.  Furthermore, if $u_1s_i> u_1$ then
\begin{equation*}
  \le(\pi_{(u_1s_i,u_2)}^{(u_1,u_2)}\ri)^{-1}\le(B_1,B_2\ri)=gx_{u_1(\al_i)}(t)\cdot\le(u_1s_i\cdot
  \Bp,u_2\cdot \Bp\ri) \qquad (t\in \C).
\end{equation*}
\item If $s_iu < u$, then there is a unique flag $B_0\in X$ such that 
$B_0\rel{s_iu} B_2$ and for all other $B'\in X$, $B'\rel{u} B_2$.  The map 
\begin{equation*}
t\mapsto gy_{u_1(\al_i)}(t)u_1\cdot \Bp
\end{equation*}
is an isomorphism $\C^*\to X-B_0$, and $B_0=gu_1s_i\cdot\Bp$.
Furthermore, if $u_1s_i> u_1$ then
\begin{equation*}
  \le(\pi_{(u_1s_i,u_2)}^{(u_1,u_2)}\ri)^{-1}(B_1,B_2)-B_0=gy_{u_1(\al_i)}(t)\cdot \le(u_1\cdot
  \Bp,u_2\cdot \Bp\ri)\qquad (t\in \C^*).
\end{equation*}
\end{enumerate}
\end{prop}

In Proposition~\ref{Bruhat-para}, the term $u_1$ can be 
incorporated into $g$, so that 
\begin{equation*}
g\cdot (u_1 \Bp,u_2 \Bp) =g u_1 \cdot (\Bp, u \Bp)
\end{equation*}
which somewhat simplifies the statement. We state the 
proposition as above since this is the form in which it will 
be used later. The proof is straightforward and is left to the reader.

\section{Double distinguished subexpressions}
\label{sec:ddses}

In this section, we will discuss the combinatorics of the ``double Weyl
group'' $W\times W$.  As with ``double flag variety,'' this terminology
can be justified by the fact that this is the Weyl group of the double
$G\times G$.  In particular, we will generalize Deodhar's notion of a
{\em distinguished subexpression}\/ to the double case.  In
Section~\ref{sec:strat-gdb-uvw}, we will see that these \ddses play the
same role for the varieties $\gdb uvw$ that distinguished subexpressions
do for $\Ri vw$.

First let us describe briefly the ``single Weyl group'' case: Consider a
reduced expression $w=s_{i_1}\cdots s_{i_n}$ in $W$, and let
$w_{(k)}=s_{i_1}\cdots s_{i_k}$, and $\Bw=\le(w_{(0)},\ldots,\wn\ri)$.
A sequence $\Bv=\le(v_{(0)},\ldots,v_{(n)}\ri)$ is called a {\bf
  subexpression} of $\Bw$ if $v_{(0)}=e$ and
\begin{equation*}
   \vk\in\{\vkm,\vkm s_{i_k}\}
\end{equation*}
for all $1\leq k\leq n$.
Informally, $\Bv$ has been obtained by throwing some of the simple
reflections out of the word $\Bw$.  

Such an expression is called {\bf distinguished} if whenever $\vkm
s_{i_k}<\vkm$, we have $\vk=\vkm s_{i_k}$, that is, $\Bv$ decreases in
length whenever possible.

Following the convention of Fomin and Zelevinsky, we consider $W\times
W$ as a Coxeter group with simple reflections $s_{-i},s_i$ for
$i\in\Pi$, the simple roots of $G$.

Pick a reduced word $s_{i_1}\cdots s_{i_n}$ for $(v,w)\in W\times W$. As
before, we let 
\begin{equation*}
\le(v_{(k)},w_{(k)}\ri)=s_{i_1}\cdots s_{i_k}.
\end{equation*}  
We let
$\ep(k)=1$ if $i_k>0$ and $-1$ if $i_k <0$.

Now, fix a sequence $\Bu=\le(u_{(0)},\ldots, u_{(n)}\ri)$ with $u_{(i)}\in
W$.  We call $\Bu$ a {\bf double subexpression of
  $(\Bv,\Bw)$} if $u_{(0)}=e$ and 
\begin{enumerate}
\item if $i_k > 0$, then $u_{(k)}\in\{u_{(k-1)},u_{(k-1)}s_{|i_k|}\}$,
\item if $i_k < 0$, then $u_{(k)}\in\{u_{(k-1)},s_{|i_k|}u_{(k-1)}\}$.
\end{enumerate}

We can write this more compactly as 
\renc{\labelitemi}{$(*)$}
\begin{itemize}
\item $u_{(k)}^{\ep(k)}\in\left\{u_{(k-1)}^{\ep(k)},u_{(k-1)}^{\ep(k)}s_{|i_k|}\right\}$.
\end{itemize}
\renc{\labelitemi}{$\bullet$}

We call $\Bu$ {\bf double distinguished} if
$u_{(k)}^{\ep(k)}=u_{(k-1)}^{\ep(k)}s_{|i_k|}$ for all $k$ such that 
$u_{(k-1)}^{\ep(k)}s_{|i_k|} < u_{(k-1)}^{\ep(k)}$. 

For each expression $\Bu$, we let
\begin{itemize}
\item $\Jo u$ be the set of indices $k\in [1,n]$ such that
  $u_{(k-1)}=u_{(k)}$,
\item $\Jp u$ be the set of indices $k\in [1,n]$ such that
  $u_{(k-1)}<u_{(k)}$,
\item $\Jm u$ be the set of indices $k\in [1,n]$ such that
  $u_{(k-1)}>u_{(k)}$. 
\end{itemize}
Obviously, a subexpression $\Bu$ of $(\Bv,\Bw)$
is uniquely determined by these subsets.

Note that each \ddse writes $\un$ as a product of $\un=\le(\uw\ri)^{-1}\uv$ with
$\uw\leq w,\uv\leq v$.


We call a \ddse \textbf{positive} if $\Jm u=\emptyset$.  As is the case
for usual (single) subexpressions, a positive subexpression with $\un=u$
exists if and only if any does.
\begin{prop}\label{sub-equiv}
    The following are equivalent:
    \renc{\labelenumi}{(\alph{enumi})}
    \begin{enumerate}
    \item There is a unique positive \ddse of any reduced expression for
    $(v,w)$ with $\un=u$.
    \item There is a \ddse of some reduced expression for
    $(v,w)$ with $\un=u$.
    \item $u=v'w'$ for $v'\leq v^{-1}$ and $w'\leq w$.
    \end{enumerate}
\end{prop}
\begin{proof}
  The implications $(a)\Rightarrow(b)\Rightarrow(c)$ are obvious.

  Now, we consider the implication $(c)\Rightarrow(a)$. 

  Since $v'< v^{-1}$ and $w'< w$, there are positive subexpressions for
  both these elements.  We can combine these into a double expression 
  $v'w'$, which is a subexpression $\bu$ for $u$.  We can
  assume by the deletion property of Coxeter groups that $\bu$ is
  positive (i.e.~reduced).  The only question is whether it is
  distinguished.
  
  Let $k$ be the largest index for which $\ukm^{\ep(k-1)} >
  \ukm^{\ep(k-1)}s_{|i_{k}|}$. Using the deletion property again, there
  is an index $k'< k$ for which $u_{(k')} > u_{(k'-1)}$, such that if we
  delete $k'$ from our original subexpression and add $k$, we have a
  positive double subexpression which agrees with our original for all
  indices outside $k'\leq i <k $ (in particular $n$), and for which the
  largest index where it is not distinguished is strictly smaller than
  $k$.  Using our cancellation argument inductively, we obtain a
  positive \ddse for $u$.

  Uniqueness follows from the fact that if a positive \ddse exists, it
  is forced by the distinguishment and positivity conditions to be the
  following: we set $\un=u$ and let
  $\tukm^{\ep(k)}=\min\le(\tuk^{\ep(k)}, \tuk^{\ep(k)}s_{|i_k|}\ri)$. 
\qed\end{proof}

It is worth noting that when either of $v=e$ (resp. $w=e$), condition (c) of
Proposition \ref{sub-equiv} simply reduces to requiring $u<w$
(resp. $u<v^{-1}$). 

We can obtain a similar reduction when $u=w_0$. 
\begin{prop}\label{doub-schub-red}
  There exists a \ddse of $(v,w)$ with $w_0$ if and only if $vw_0 \leq w$.
\end{prop}
\begin{proof}
  By Proposition \ref{sub-equiv}, if such a \ddse exists, then
  $w_0=v'w'$ with $v'\leq v,w'\leq w$.  Thus $vw_0\leq
  (v')^{-1}w_0=w'\leq w$.

  On the other hand, if $vw_0\leq w$, then $v'=v, w'=vw_0$ realizes
  condition $(c)$ of Proposition \ref{sub-equiv}.\qed
\end{proof}

\section{The stratification of $\gdb uvw$}
\label{sec:strat-gdb-uvw}

For each reduced expression $(\Bv,\Bw)$ for $(v,w)$, we have the standard
projection maps
\begin{equation*}
\pi_k: \flag_v\times\flag^w\to \flag_{\vk}\times\flag^{\wk}, \quad
\pi_k = \pi^v_{\vk} \times \pi^w_{\wk}.
\end{equation*}

Fix such an expression, and for a pair $(B_1,B_2)\in\gdb
uvw\subset\flag_v\times\flag^{w}$, consider the sequence
$\uk=r(\pi_k(B_1,B_2))w_0$.  That is, $\uk$ is the unique element such
that $\pi_k(B_1,B_2)\in\gdb \uk\vk\wk$.

\begin{thm}
  The sequence $\Bu=(u_{(0)},\ldots,\un)$ is a \ddse of $(\Bv,\Bw)$.
\end{thm}
\begin{proof}
To simplify notation, we let
\begin{equation*}
\pi_k^1=\pi^v_{\vk} :  \flag_v \to \flag_{\vk}, \quad
\pi_k^2=\pi^w_{\wk} :  \flag_w \to \flag_{\wk}.
\end{equation*}

  Assume for simplicity that $i_k<0$ (the proof for $i_k>0$ is the
  same).  Then $\pi_k^2(B_2)=\pi_{k-1}^2(B_2)$, and
  $r(\pi_k^1(B_1),\pi_{k-1}^1(B_1))=s_{|i_k|}$. 

  By Proposition~\ref{Bruhat-para}, $\uk\in\{\ukm,s_{|i_k|}\ukm\}$, so
  $\Bu$ is a double subexpression of $(\Bv,\Bw)$.  If  $\ukm >
  s_{|i_k|}\ukm$, since left multiplication by $w_0$ is order-reversing
  for Bruhat order, 
  \begin{equation*}
    r\le(\pi_k^1(B_1),\pi_{k-1}^1(B_1)\ri)\cdot r\le(\pi_{k-1}(B_1,B_2)\ri)=r\le(\pi_k(B_1,B_2)\ri)
  \end{equation*}
  by property (\ref{item:4}) of Proposition~\ref{prop:rel-pos-prop}.
  That is, $\uk=s_{|i_k|}\ukm$.  Thus, $\Bu$ is also distinguished.
\qed\end{proof}

Let $\gdbb uvw$ to be the subvariety of $\gdb uvw$ defined by  
\begin{equation*}
  \gdbb uvw = \{(B_1,B_2)\vert \pi_k(B_1,B_2)\in \gdb {\uk}{\vk}{\wk}\}
\end{equation*}
The theorem above shows that $\gdb uvw$ has a partition, 
\begin{equation}\label{eq:3}
  \gdb uvw=\bigsqcup_{\Bu}\gdbb uvw
\end{equation}
which we will show in Section~\ref{sec:coisotr-deodh-strat} is a
stratification. First, let us understand the topology of the strata a
bit better.

We let 
  \begin{equation*}
    g_k=
    \begin{cases}
       y_{\uk^{\Bv}\al_{i_k}}(t_k)& i_k>0, k\notin\Jp u\\
      x_{\uk^{\Bw}\al_{i_k}}(t_k)& i_k< 0, k\notin\Jp u\\
      1 & k\in\Jp u
    \end{cases}
  \end{equation*}
and let $g_{\Bu}(t_1,\ldots, t_n)=g_1(t_1)\cdots g_n(t_n)$ where
$t_i=1$ if $i\in\Jp u$, $t_i\in\C^*$ if $i\in\Jo u$ and $t_i\in \C$ if
$i\in \Jm u$.

\begin{thm}\label{thm:strata-param}
  For any \ddse $\Bu$, 
  \begin{equation*}
    \gdbb uvw\cong (\C^*)^{|\Jo u|}\times \C^{|\Jm u|}.
  \end{equation*}
  In particular, $\dim \gdbb uvw = n-|\Jp u|\leq\ell(w)+\ell(v)-\ell(u)$.

  Furthermore, the map 
  \begin{equation*}
  \vp_{\Bu}:(t_1,\ldots,t_n)\mapsto g_u(t_1,\ldots,t_n)\cdot 
  (\uw \cdot \Bp ,\uv\cdot \Bm )
  \end{equation*}
  is an isomorphism between $\CsC{|\Jo u|}{|\Jm u|}$ and $\gdbb uvw$.
\end{thm}
\begin{cor}\label{cor:nonempty}
  The set $\gdb uvw$ is a smooth, locally closed, and irreducible
subvariety of $\flag \times \flag$. It is nonempty if and only
if $u = w' v'$ for some $w' \leq w^{-1}$ and $v' \leq v$.

\end{cor}
The first part of this corollary is a special 
case of \cite[Proposition 4.2]{EL04}. That part and 
Proposition 4.2 in \cite{EL04} also follow directly 
from Richardson's result \cite[Corollary 1.5]{Ri92}. 
\begin{proof}
  The intersection between $\Bp\times\Bm$ and $G_\De$-orbits is
  transverse, since $\fr g\oplus\fr g$ is spanned by $\fr b_+\oplus\fr
  b_-$ and $\fr g_\De$.  Thus, $\gdb uvw$ is a locally closed,  
  smooth subvariety of $\flag \times \flag$.
  of dimension $\ell(v)+\ell(w)-\ell(u)$.  In particular, in each irreducible 
  component,
  there is a stratum of dimension $\ell(v)+\ell(w)-\ell(u)$.  Since
  \begin{equation*}
    \dim\gdbb uvw=\ell(v,w)-|\Jp u|=\ell(v,w)-\ell(u)-|\Jm u|
  \end{equation*}
  a stratum is of maximal dimension if $\Bu$ is positive, and there is a
  unique such \ddse.  Thus, $\gdb uvw$ only has one component, and thus
  is irreducible.

  The second part follows from Proposition~\ref{sub-equiv} 
and Theorem~\ref{thm:strata-param}. 
\qed\end{proof}

Applying Proposition~\ref{double-schub}, we see that
Corollary~\ref{cor:nonempty} also gives a characterization of when the
varieties $\Ri vw$ are nonempty (in fact, it provides 3 different
descriptions).  Proposition~\ref{doub-schub-red} and the discussion
preceding it show that these characterizations all reduce to a
previously known characterization.
\begin{cor} \emph{(\cite[Corollary~1.2]{Deo85})}
  The variety $\Ri vw$ is nonempty if and only if $v\leq w$.
\end{cor}

\begin{proof}[Proof of Theorem~\ref{thm:strata-param}]
  By induction, assume the result is true for $\unm,\vnm,\wnm$, and for
  simplicity, assume $i_n>0$, and let 
  \begin{equation*}
    X=\{(B_1,B_2)\in \flag\times\flag |
  B_1\rel{s_{i_n}} B_1', (B_1',B_2)\in\gdbb \unm\vnm\wnm\}. 
  \end{equation*}
Then, if $s_{i_n}\unm< \unm$, then $X=\gdbb uvw$, and by
Proposition~\ref{Bruhat-para}, part (1), 
\begin{multline*}
  g_u(t_1,\ldots, t_n)\cdot (\uw \cdot \Bp ,\uv\cdot \Bm
  )\\
=g_{\unm}(t_1,\ldots, t_{n-1})\cdot y_{\uw(\al_{i_n})}(t_n)\cdot (s_i\uw_{(n-1)} \cdot \Bp ,\uv\cdot \Bm)
\end{multline*}
is a parameterization of $\gdbb uvw$.

On the other hand, if $s_{i_n}\unm> \unm$, then $X=\gdbb \unm v w\cup
\gdbb {s_{\mathit{i_n}}\unm} v w$, and by a similar calculation, part
(2) of Proposition~\ref{Bruhat-para} confirms that $g_u$ provides a
parameterization.
\qed\end{proof}

If $u=e$, then the unique non-decreasing sequence is obviously given
by $e^+_{(k)}=e$ for all $k$.

In \cite{FZ99}, Fomin and Zelevinsky construct a dense subset of
$G_{v,w}$ for each reduced word for $(v,w)$ by the factorization map 
\begin{align*}
\xvw:H\times(\C^*)^{n}&\to G_{v,w}\\
  (h,t_1,\ldots t_n)&\mapsto x_{i_1}(t_1)\cdots x_{i_n}(t_n)h.
\end{align*}

Clearly, the composition $\vp\circ \xvw$ does not depend on $h$, and so
gives an injection $(\C^*)^n\to \gdb {e}vw$.

\begin{cor}
  $\gdbb {e^+}vw$ is exactly the set of elements of the
  form $(g\cdot \Bp ,g\cdot \Bm )$ where $g=h x_{i_1}(t_1)\cdots
  x_{i_n}(t_n)$ with $t_i\in\C^*$.
\end{cor}

Thus the stratification of $\gdb {e}vw$ induces a stratification 
of the double Bruhat cell $\Bp v \Bp \cap \Bm w \Bm$ with maximal 
stratum that coincide with the open subsets of Fomin-Zelevinsky
\cite{FZ99}. In addition, all strata have dimension greater than more
half the dimension of $G^{u,v}$; 
in fact, they are coisotropic with respect to the standard
Poisson structure on $G$, as we will show in the 
last section.

Furthermore, if $w=1$, then we recover a rewriting of the
parameterizations of \cite{MR04} for double Schubert cells: to produce
our parameterization from theirs, simply commute all Weyl group elements
past the unipotent elements until they are collected at the end of the word.

On the other hand, if $v=1$, then we recover the double Schubert cell
$\mc R_{uw_0,vw_0}$,  with the opposite parameterization.  

\begin{remark}
  These results all have a straightforward generalization to the
  intersection of a $G_\De$-orbit with an orbit of $B\times B'$ with $B$
  and $B'$ two Borel subgroups of $G$, which may not be opposite. 
  Unlike the varieties
  $\gdb uvw$, we do not expect these more general orbit intersections to
  have good properties with respect to Poisson or positive structures.
  However, in the case where $B=B'$, we will recover results of Curtis
  \cite{Cur88}, which are relevant to the study of Hecke algebras.

  The varieties of Curtis appear in a second way in our theory.  If one
  chooses a reduced word for $(v,w)$ which puts all reflections
  in the first copy of $W$ before those appear in the second copy,
  i.e. an ``unmixed'' word, then $\pi_{\ell(v)}$ is simply the
  projection onto the first factor $\gdb uvw \subset
  \flag\times\flag\to\flag$.  Each fiber of this map will be isomorphic
  to one of Curtis's varieties, though which variety it is will depend
  on which double Bruhat cell we are taking the fiber over.  In this
  case, each of our strata projects to a Deodhar stratum in $\flag$,
  with fiber given by one of Curtis's strata.  In particular, considering
  the intersection of our stratification with each fiber will recover
  Curtis's stratification.
\end{remark}

\section{Chamber minors}
\label{sec:chamber-minors}

Just as in \cite{BZ97,MR04}, we can define $\gdbb uvw$ in terms of
certain generalized minors, which one can call ``generalized chamber
minors.'' Although Marsh and Rietsch have already claimed this name for
the generalized minors which appear in their ``generalized Chamber
Ansatz,'' since our situation subsumes theirs, it should cause no
confusion.  

For a fixed highest weight
vector $h_\la$ in the representation $V_\la$, let $h_\la^*$ denote the
unique dual \hw vector in $V_\la^*$, and $\inn{-,-}$ denote the standard
pairing between these spaces.  The generalized minors
$\Delta^{w\la}_{v\la}$ are defined by
\begin{equation*}
  \Delta^{w\la}_{v\la}(g)=\inn{g\bar{v}\cdot h_\la,\bar{w}\cdot h_\la^*}.
\end{equation*}
These are precisely the matrix coefficients of extremal weight vectors
in representations of $G$.  They play a central role in the papers
\cite{BFZ05,BZ97,FZ99,KZ02,Zel00}.

Let $\Up$ and $\Um$ be the unipotent radicals of $\Bp$ and $\Bm$, and
for all $w\in W$, let $U^{\pm}_w=U^{\pm}\cap wU^{\mp}w^{-1}$.  
As usual, the map $\al_{v,w}:\Up_v\times\Um_w\to \flag_v\times\flag^w$ given by
action on $(v\cdot \Bp,w\cdot\Bm)$ is an isomorphism.  Furthermore, 
\begin{prop}
  For all $B_1\in\flag_w,B_2\in \flag^v$, 
  \begin{equation*}
    \pi_k(B_1,B_2)=\al^{-1}_{v,w}(B_1,B_2)\cdot (\vk\cdot \Bp,\wk\cdot\Bm).
  \end{equation*}
\end{prop}

This allows us to identify the varieties $\gdbb uvw$ using chamber
minors:
\begin{thm}
  Let $(z_1,z_2)=\al^{-1}_{v,w}(B_1,B_2)$.  Then $(B_1,B_2)\in \gdbb uvw$ if and
  only if 
  \begin{align*}
    \gm {\vk\uk w_0}{\wk w_0}{i_k}(z_1^{-1}z_2)\neq 0 \text{ for all }
    k\in\Jo u, i_k<0,\\
    \gm {\vk\ukm w_0}{\wk w_0}{i_k}(z_1^{-1}z_2) = 0 \text{ for all } k\in\Jp
    u, i_k <0,\\
    \gm {\wk\uk^{-1}}{\vk}{i_k}(z_2^{-1}z_1)\neq 0 \text{ for all }
    k\in\Jo u, i_k>0,\\
    \gm {\wk\ukm^{-1}}{\vk}{i_k}(z_2^{-1}z_1)= 0 \text{ for all }
    k\in\Jp u, i_k>0.\\
  \end{align*}
\end{thm}
\begin{proof}
  By definition $(B_1,B_2)=(z_1v\cdot\Bp,z_2w\cdot\Bm)$ and
  \begin{equation*}\pi_k(B_1,B_2)=(z_1\vk\cdot\Bp,z_2\wk\cdot\Bm),\end{equation*} so for all $k$, we calculate  
  \begin{equation*}
    r(\pi_k(B_1,B_2))=r(z_1\vk\cdot\Bp,z_2\wk\cdot\Bm)=r(\Bp,
    \vk^{-1}z_1^{-1}z_2\wk\cdot \Bm)
  \end{equation*}
  is $\uk w_0$ if and only if $\vk^{-1}z_1^{-1}z_2\wk\in \Bp \uk \Bm$.

  By basic representation theory, $g\in \Bp u \Bm$ if and only if 
  \begin{align*}
    \gm {uw_0}{w_0}i(g)&\neq 0\\
    \gm {u'w_0}{w_0}i(g)&
      =0 \qquad u'\om_i < u\om_i
  \end{align*}
  or equivalently
 \begin{align*}
    \gm {u^{-1}}{}i(g^{-1})&\neq 0\\
    \gm {u'}{}i(g^{-1})&
      =0 \qquad u'\om_i < u^{-1}\om_i.
  \end{align*}
  Applying these to $g= \vk^{-1}z_1^{-1}z_2\wk$, we find that the
  relations are necessary.

  Now assume that they are sufficient for all $(v',w')$ with
  $\ell(v',w')<\ell(v,w)$. Thus, the first $n-1$ relations are
  sufficient to assure that
  $\pi_{n-1}(B_1,B_2)\in\gdbb\unm\vnm\wnm$. For ease, again assume that
  $i_n>0$.

  If $n\in \Jm u$, we are done, so assume $n\in\Jo u\cup\Jp u$, that is
  $s_{|i_n|}\unm> \unm$. Then $(B_1,B_2)$ must be in $\gdb \unm vw$ or
  $\gdb {s_{|i_n|}\unm }vw$, and this determines which Deodhar stratum
  it lives in.   Since $\unm^{-1}s_{|i_n|}\om_i>\unm^{-1}\om_i$, if
  $\un=s_{|i_n|}\unm$, then by the formulae above, 
  $\gm {w\unm^{-1}}{v}i(z_2^{-1}z_1)=0$, whereas if $\un=\unm$ then 
  \begin{equation*}
    \gm {w\unm^{-1}}{v}{i}(z_2^{-1}z_1)=\gm
    {wu^{-1}}{v}{i}(z_2^{-1}z_1)\neq 0.
  \end{equation*}
  Thus, the hypotheses of the theorem are sufficient as well as necessary.
\qed\end{proof}

In the case where $w=1$, there is an explicit change of coordinates
between our parameterization of Theorem~\ref{thm:strata-param}, and a
coordinate system given by chamber minors, known as the Generalized
Chamber Ansatz, described in \cite{MR04}.

\begin{question}
Is there a generalization of the Chamber Ansatz to the coordinate systems in Theorem~\ref{thm:strata-param}?
\end{question}

\section{Positivity}
\label{sec:positivity-1}

Previous work along these lines has been closely related to the theory
of total positivity.  While this the concept of totally positive
matrices has existed for decades, it was developed in its modern form by
Lusztig, followed by Berenstein, Fomin and Zelevinsky (see the papers
\cite{BZ97,BZ01,FZ99,Lus94,Lus98a}). 

Let $\Gsp$, the strictly totally positive part of $G$, be the subset of $G$ on
which all generalized minors are positive.  By \cite{FZ99}, this is the
same as the set of $G$ of the form
\begin{equation*}
  g=x_{i_1}(t_1)\cdots x_{i_n}(t_n)
\end{equation*}
for $t_i\in \R_{>0}$ and $s_{i_1}\cdots s_{i_n}$ a reduced word for $(w_0,w_0)$.  

The set of
non-negative elements of $G$ is simply the closure
$\Gp=\overline{\Gsp}$. Similarly, in any based $G$-space 
$(X,x)$, we can define a positive subset relative to $x$ as the subset $\pos{X}=\overline{\Gsp\cdot
  x}$.

In \cite{Lus94}, Lustig showed that $\Gp$ has a natural
(real) cell decomposition, in which each cell is the intersection of
$\Gp$ with a double Bruhat cell.

In \cite{R99}, Rietsch showed that $\pos\flag$ has a similar
cell decomposition in terms of double Schubert cells, as was conjectured
by Lusztig.

Now, we define {\em the positive part of $\flag\times \flag$} to be 
$\pos{(\flag\times\flag)} = \overline{G_\De^{+} \cdot (\Bp, \Bm)}$. We prove 
below that this is not the product of the positive flag varieties
$\pos\flag\times\pos\flag\subset \flag\times\flag$. In particular,
it differs from the restriction of the nonnegative part of the 
wonderful compactification $\overline{G}$, 
defined by Lusztig \cite{Lus98b}, intersected with the lowest stratum 
$\flag \times \flag$. In the terminology
of \cite{FG03b}, these different ``positive parts'' correspond to
different positive structures on $\flag\times\flag$, one the product the
standard positive structures on the flag variety, and one descending from the
standard positive structure on $G$.   

For example, if
$G=\SL 2\C$, one can identify $\flag\times\flag$ with
$\mathbb{CP}^1\times\mathbb{CP}^1$. Then $\pos\flag\times\pos\flag$ is the
subset 
\begin{equation*}
  \{(a,b)|a,b\in[0,\infty]\}
\end{equation*}
whereas $\pos{(\flag\times\flag)}$ is 
\begin{equation*}
  \{(a,b)|a\in[0,\infty],b\in [1/a,\infty]\}.
\end{equation*}
While this definition may look asymmetric in $a$ and $b$, in fact, it is
invariant under the switch map.

\begin{thm}
  For any simple algebraic group $G$,
  \begin{equation*}
    \pos{(\flag\times\flag)}\subsetneq \pos{\flag}\times\pos\flag.
  \end{equation*}
\end{thm}
\begin{proof}
  Since $\Gsp=\pos{\Up}\pos{H}\pos{\Um}=\pos{\Um}\pos{H}\pos{\Up}$ by
  \cite{FZ99}, 
  \begin{equation*}
    \Gsp\cdot(\Bp,\Bm)\subset \pos{\Um}\cdot\Bp\times\pos{\Up}\cdot\Bm.
  \end{equation*}
  Taking closure of both sides yields the desired inclusion.  
  
  Now, consider the Schubert cell
  $\overline{\flag_{s_i}\times\flag^{s_i}}$ for some simple root $s_i$.
  This naturally identified with the flag variety of the corresponding
  root subgroup, isomorphic to $\SL 2\C$. The intersection
  $\pos{\flag}\times\pos\flag\cap
  \overline{\flag_{s_i}\times\flag^{s_i}}$ is precisely the same as
  the corresponding product in the double flag of $\SL 2\C$ (this is clear
  from the parameterization of $\flag$ given in \cite{MR04}).  On the
  other hand $\pos{(\flag\times\flag)}\cap\gdb{e}{s_i}{s_i}$ is
  precisely $\pos{(G^{s_i,s_i})}\cdot(\Bp,\Bm)$.  By our calculation for
  $\SL 2\C$, these sets do not coincide.
\qed\end{proof}

By analogy with previous positivity results \cite{FZ99, He04, MR04}, 
we conjecture that the
varieties $\gdb uvw$ and the parameterizations of them we have described
are closely connected to the positive part of $\flag\times\flag$. 

\begin{conj}
  $\pos{(\flag\times\flag)}$ has a real cell decomposition, given by its
  intersections with the subvarieties $\gdb uvw$ where $u,v,w$ range
  over $W$.  Each such intersection $\pos{\gdb uvw}=\gdb
  uvw\cap\pos{(\flag\times\flag)}$ is contained in $\gdbb {u^+}vw$, and
  $\vp_{\Bu^+}$ gives an isomorphism between
  $\R_{>0}^{\ell(w)+\ell(v)-\ell(w)}$ and $\pos{\gdb uvw}$.
\end{conj}

The positive structure on $G$ can be strengthened to a cluster algebra
structure (for background on cluster algebras, see
\cite{BFZ05,FG03b,FZ02}), so the varieties $\gdb evw=L^{v,w}$ have a
natural cluster algebra structure.  
\begin{question}\label{q:clusters}
Does $\gdbb uvw$ have a cluster algebra structure for $u\neq e$? 
\end{question}
While quite natural from the perspective that the varieties $\gdbb uvw$
are a generalization of double Bruhat cells, this question appears to be
quite difficult, and is beyond the scope of this paper.

\section{The Poisson geometry of strata}
\label{sec:coisotr-deodh-strat}

In this section, we show how our partition (and that of Deodhar) are
compatible with the natural Poisson structure induced by the standard
Poisson structure on $G$ and the double structure on $G\times G$.  We
assume from now on that $G$ and $G\times G$ are endowed with these
Poisson structures, unless stated otherwise.

There seems to be a clear, if somewhat imprecise, connection between
positive and Poisson structures on $\flag\times\flag$.  Just as
$\flag\times\flag$ has two positive structures, 
\begin{itemize}
\item one coming from the product of standard positive structures on $\flag$ 
\item one from the standard positive structure on $G_\De$, 
\end{itemize}
it also has two
Poisson structures, 
\begin{itemize}
\item one from the product of the Poisson structure induced on
  $\flag$ by the standard Poisson group structure on $G$ (``the product
  structure''). 
\item one induced by the double Poisson structure on $G\times G$ (``the
  double structure'').
\end{itemize}

One possible explanation for this correlation would be a positive answer to Question~\ref{q:clusters}.  A cluster algebra
structure on the varieties $\gdb uvw$ would endow them
with compatible Poisson and positive structures. 

In fact, there are many different Poisson structures on
$\flag\times\flag$, in particular, one for each pair of transverse
Lagrangian subalgebras $\fr u$ and $\fr u'$
of $\fr g\times\fr g$.  For more on such
subalgebras, see the papers of Evens and Lu, \cite{EL01,EL04}.
It is proved in \cite{LY06} that the intersection of each
$N(\fr u)$ and $N(\fr u')$ orbits is a regular Poisson subvariety
of $\flag \times \flag$, where $N(\fr u)$ and $N(\fr u')$
are the normalizers of ${\fr u}$ and ${\fr u'}$ 
in $G \times G$. This provides the analog of the partition 
\eqref{0:part} of $\flag \times \flag$ in the general case.
\begin{question}
Is there a generalization of the 
finer partition \eqref{eq:3} to general choices of transverse Lagrangian subalgebras which has good Poisson properties? 
\end{question}

Let us first collect a few useful pieces of information about
Poisson algebraic groups.
\begin{prop}
  The subgroups $\Bp,\Bm$ are Poisson subgroups of $G$, and $\Bp\times\Bm$ and $G_\De$ are
  Poisson subgroups of $G\times G$.
\end{prop}

This implies that the standard Poisson structure on $G$ and 
the double structure on $G\times G$ can be pushed forward to well
defined Poisson structures on $\flag$ and
$\flag\times\flag$. The latter become Poisson homogeneous spaces 
for $\flag$ and $\flag\times\flag$, respectively, of so called group type.

\begin{prop}
  Let $X$ be a Poisson $H$-variety for any Poisson algebraic group $H$, and
  assume the Poisson tensor vanishes at $x\in X$.  Then the map $h\mapsto
  h\cdot x$ is Poisson.  In particular, $H\cdot x$ is a Poisson
  subvariety.
\end{prop}
\begin{proof}
  The Poisson tensor vanishes at $x$, so $g\mapsto
  (g,x)$ is a Poisson map. The rest is clear. Here again
$H \cdot x$ becomes a Poisson homogeneous $H$-space of group type. 
\qed\end{proof}

This proposition implies that $\Ri vw$ and $\gdb uvw$ are Poisson
subvarieties of $\flag$ and $\flag\times\flag$.  In fact, they are
torus orbits of symplectic leaves (see \cite{EL04,GY05}).

It also implies that many projection maps between orbits are Poisson.  
\begin{prop}
  The projection maps 
  \begin{align*}
    \pi^w_{w'}&:\flag_w\to\flag_{w'}\\
    \bar {\pi}^w_{w'}&:\flag^w\to\flag^{w'}
  \end{align*}
  for $w,w'\in W$ with $\ell(w^{-1}w')=\ell(w)-\ell(w')$ are Poisson submersions, 
  as are the analogous maps for $\Bp\times\Bm$ orbits on $\flag\times\flag$.
\end{prop}
Since the Poisson structure on $\flag\times\flag$ is not the product of
Poisson structures on $\flag$, the second half of the proposition does
not follow immediately from the first, though it does have the same
proof. Similar techniques were used in \cite[Proposition 1.6, \S 3.1]{GY05}.
\begin{proof}
  If $f:X\to Y$ and $g:Y\to Z$ are surjective maps between Poisson
  varieties, and $f$ and $g\circ f$ are both Poisson then $g$ is Poisson
  as well.  

  Now, let $f$ be the projection map $\Bp\to \Bp w\cdot \Bp$
  given by action on $w\cdot \Bp$.  Since the Poisson tensor vanishes at
  $w\cdot\Bp$, this map is Poisson.

  Let $g=\pi^w_{w'}$.  In this case, $g\circ f$ is just the map
  $\Bp\to\Bp w'\cdot \Bp$ given by action on $w'\cdot \Bp$, which as we
  argued above, is Poisson. Thus, $\pi^w_{w'}$ is Poisson as well.

  The same argument works on $\flag\times\flag$.
\qed\end{proof}

We might hope that there is some sort of compatibility between the
Poisson structure of $\gdb uvw$ and our subvarieties $\gdbb uvw$.  They
cannot be Poisson since $\gdb uvw$ are the minimal $T$-invariant Poisson
subvarieties of $\flag\times\flag$.  But there is a weaker notion of
compatibility with Poisson structures which it is more reasonable to expect. 

One calls a subscheme $Y$ of an algebraic Poisson variety $X$ 
{\bf coisotropic} if its ideal sheaf $\mc I_Y$ is closed under Poisson
bracket.  That is, if $f$ and $g$ are two rational functions on $X$,
which both vanish on $Y$, $\{f,g\}$ does as well.  Note that this is a
much weaker condition than being a Poisson subvariety, which requires
$\{f,g\}$ to vanish on $Y$ if either of $f$ or $g$ does.

\begin{prop}\label{prop:coisotropy}
  If $Y\subset X$ is coisotropic and $\psi:Z\to X$ is a Poisson map,
  $\psi^{-1}(Y)$ is coisotropic in $Z$.
\end{prop}

Some care must be taken when applying this proposition to algebraic maps
which are not reduced.  In this case, one may pullback the reduced
subscheme structure on a subset of $X$, and obtain a non-reduced
subscheme $\psi^{-1}(Y)$ which is coisotropic, even though its reduced
counterpart $\psi^{-1}(Y)_{red}\subset Z$ is not.  This will not be a
problem as long as the map $\psi$ is smooth (in the algebraic sense,
which is roughly equivalent to being a submersion in the differential
category).

\begin{proof}
  Obviously, this only needs to be checked for affine varieties, so let
  $\vp^*:A\to B$ be a homomorphism of Poisson algebras, and $I\subset A$
  be an ideal such that $\{I,I\}\subset I$, then $\{\vp(I),\vp(I)\}
  \subset \vp(I)$, and so 
  \begin{equation*}
    \{B\vp(I),B \vp(I)\}\subset B
  \{\vp(I),\vp(I)\}+B\vp(I)\{\vp(I),B\}+\vp(I)\cdot \{B,B\}\subset B\vp(I)\qed
  \end{equation*}
\end{proof}

\begin{thm}
  The Deodhar components $\Rb vw$ and $\gdbb uvw$ are coisotropic
  subvarieties of $\flag$ and $\flag\times\flag$ respectively.
\end{thm}
\begin{proof}
  We will first consider the case of $\Rb vw$. By
  Proposition~\ref{prop:coisotropy}, coisotropy is preserved under
  pullback by Poisson maps. Thus, consider the Poisson map
  $\pi_{n-1}:\flag_{\wn}\to \flag_{\wnm}$.  Each Deodhar component of
  $\Ri vw$ is the intersection of the pull-back of a Deodhar component
  of $\Ri {vs_{i_n}}{\wnm}$ or $\Ri {v}{\wnm}$ with $\flag^{w_0v}$ (here
  we use the fact that $\pi_{n-1}$ is smooth (in the sense of algebraic
  geometry), to see that the subscheme structures match
  up).  By induction, these are coisotropic subvarieties of
  $\flag_{\wnm}$, and so their pullbacks are coisotropic in
  $\flag_{\wn}$, which is a Poisson subvariety of $\flag$.

  The proof for $\gdb uvw$ is precisely the same.
\qed\end{proof}

We complete this section with a proof of the fact that our partition is
in fact a stratification.  While the dimension calculations we have done
suggest that this is the case, it needs to be carefully checked.  In
some other, superficially similar situations, e.g. certain moduli spaces
of flat connections, similar partitions are not stratifications.
Luckily, on $\flag\times\flag$, this is not the case.

\begin{thm}
  The partition in equation~\eqref{eq:3} is a stratification of
  $\flag_v\times\flag^w$ (and thus of $\gdb uvw$). 
\end{thm}
\begin{proof}
  Let $\mc X$ be the sheaf of all algebraic vector fields on
  $\flag_v\times\flag^w$ and let $\mathcal{G}_i$ the subsheaf of vector
  fields such $X$ such that the pushforward $T\pi_k(X)$ is well defined
  and Hamiltonian on $\flag_{\vk}\times\flag^{\wk}$.  In particular,
  $\mc G_n$ is simply the sheaf of Hamiltonian vector fields.  Now let
  $\mc G$ be the intersection of these sheaves.  This sheaf is $H$
  invariant, the maps $\pi_k$ and the Poisson structure on each variety
  are $H$-invariant, so $[\fr h,\mc G]\subset \mc G$, and $\mc G'=\fr
  h+\mc G$ is a Lie algebra subsheaf of $\mc X$. Consider the orbits of
  $\mc G'$ on $\flag_v\times\flag^w$, that is, the foliation one obtains
  from integrating it.  This is a partition of $\flag_v\times\flag^w$
  which respects the partition $\flag_v\times\flag^w=\sqcup\gdb uvw$
  (since these are torus orbits of symplectic leaves).  We claim that
  these are precisely the subvarieties $\gdbb uvw$.  Since the varieties
  $\gdbb uvw$ are connected, we need only check
  that the image of the evaluation map $\mathrm{ev}_x:\mc G'\to
  T_x(\flag_v\times\flag^w)$ has image precisely $T_x\gdbb uvw$.  That the
  image is contained in this tangent space is clear, since
  $T_x\pi_k\circ \mathrm{ev}_x$ lands in the tangent space to the symplectic
  leaf in each case, since all push forwards are Hamiltonian.  On the
  other hand, assume that the claim is true for $k<n$.  Then if
  $u>\unm$, 
  \begin{equation*}
   T_x\pi_{n-1}:T_x\gdbb uvw\to T_x\gdbb {\unm}{\vnm}{\wnm}
  \end{equation*}
  is an isomorphism, and the hamiltonian vector fields $X_{\pi_k^*f}$
  for all functions $f$ on $\flag_{\vnm}\times\flag^{\wnm}$ span the
  tangent space to $\gdbb uvw$.  Otherwise, there is a 1 dimensional
  kernel, which the Hamiltonian vector fields of pullbacks form a
  complement to in $T_x\gdbb uvw$.  Since evaluation by the Poisson form
  $\Pi_{\sharp}$ vanishing on
  $\ker\, T_x\pi_{n-1}$ is an open condition, on some (analytic)
  neighborhood of $x$,  $\Pi_\sharp$ does not vanish on $\ker
  T_x\pi_{n-1}$, and there is a 1-form $\si$ such that
  $\Pi_{\sharp}\si\in\ker \,T_x\pi_{n-1}$.  The Hamiltonian vector field
  $X_\si$ is in $\mc G$, since it is killed by $T_x \pi_{k}$ for all
  $k<n$.

  Thus, we have proved that the subvarieties $\gdbb uvw$ are precisely
  the leaves of the foliation corresponding to $\mc G'$.  Since the
  closure $\overline{\gdbb uvw}$  is also invariant under $\mc G'$, it
  is a union of leaves of the corresponding foliation, i.e. of
  varieties $\gdbb {u'}vw$ where $u'$ ranges over some set of reduced
  double subexpressions of $(v,w)$.  Since the varieties  $\gdbb uvw$
  are smooth, connected and locally closed, this partition is, in fact,
  a stratification.
\qed\end{proof}

Consider the poset structure on \ddses given by $\mathbf{u'}\geq
\mathbf{u}$ if and only if $\uk'\geq\uk$ for all $k$.  Since Bruhat
order gives the closure relations on $\flag\times\flag$, the subvariety
\begin{equation*}
Q^{\mathbf{u}}_{\mathbf{v,w}}=\bigsqcup_{\mathbf{u'}\geq \mathbf{u}}\gdbb {u'}vw
\end{equation*}
is closed.  By definition, $\overline{\gdbb uvw}\subset
Q^{\mathbf{u}}_{\mathbf{v,w}}$.

\begin{question}
  Are these sets equal?  If not, what are the closure relations for the
  stratification of $\flag_v\times\flag^w$ by $\gdbb uvw$?
\end{question}

\def\cprime{$'$}

\end{document}